\documentclass[10pt]{article}
\usepackage{mathrsfs}
\usepackage{amssymb}
\usepackage{amsmath}
\usepackage[all]{xy}

\def\id{\mathop{\rm id}\nolimits}
\def\Im{\mathop{\rm Im}\nolimits}
\def\Ker{\mathop{\rm Ker}\nolimits}

\def\mod{\mathop{\rm mod}\nolimits}
\def\Mod{\mathop{\rm Mod}\nolimits}
\def\Hom{\mathop{\rm Hom}\nolimits}

\def\Ext{\mathop{\rm Ext}\nolimits}

\def\pd{\mathop{\rm pd}\nolimits}

\def\resdim{\mathop{\rm res.dim}\nolimits}
\def\coresdim{\mathop{\rm cores.dim}\nolimits}

\def\proj{\mathop{\rm proj}\nolimits}
\def\inj{\mathop{\rm inj}\nolimits}
\def\GProj{\mathop{\rm GProj}\nolimits}
\def\GInj{\mathop{\rm GInj}\nolimits}
\def\Gproj{\mathop{\rm Gproj}\nolimits}
\def\Ginj{\mathop{\rm Ginj}\nolimits}

\def\sup{\mathop{\rm sup}\nolimits}

\def\Con{\mathop{\rm Con}\nolimits}
\title{\large \bf Applications of Balanced Pairs
}
\author{Huanhuan Li\thanks{{\it E-mail address}: lihuanhuan0416@163.com},
\  Junfu Wang\thanks{{\it E-mail address}: wangjunfu05@126.com}
\ and Zhaoyong Huang\thanks{{\it E-mail address}: huangzy@nju.edu.cn}\\
{\footnotesize \it Department of Mathematics, Nanjing University,
Nanjing 210093, Jiangsu Province, China}}
\date{}
\begin{document}
\baselineskip=15pt
\maketitle
\begin{abstract}
Let $(\mathscr{X}$, $\mathscr{Y})$ be a balanced pair in an abelian category. We first introduce the notion of
cotorsion pairs relative to $(\mathscr{X}$, $\mathscr{Y})$, and then give some equivalent characterizations when
a relative cotorsion pair is hereditary or perfect. We prove that if the $\mathscr{X}$-resolution dimension
of $\mathscr{Y}$ (resp. $\mathscr{Y}$-coresolution dimension
of $\mathscr{X}$) is finite, then the bounded homotopy category of $\mathscr{Y}$ (resp. $\mathscr{X}$)
is contained in that of $\mathscr{X}$ (resp. $\mathscr{Y}$). As a consequence, we get that the right
$\mathscr{X}$-singularity category coincides with the left $\mathscr{Y}$-singularity category if
the $\mathscr{X}$-resolution dimension of $\mathscr{Y}$ and the $\mathscr{Y}$-coresolution dimension
of $\mathscr{X}$ are finite.
\end{abstract}

\vspace{0.5cm}

\centerline{\bf  1. Introduction}

\vspace{0.2cm}

Cartan and Eilenberg introduced in [CE] the notions of right and left balanced functors.
Then Enochs and Jenda generalized in [EJ1] them to relative homological algebra as follows. Let $\mathscr{C}$,
$\mathscr{D}$ and $\mathscr{E}$ be abelian categories and $T(-,-):\mathscr{C}\times\mathscr{D}\to\mathscr{E}$
be an additive functor contravariant in the first variable and covariant in the second. Then $T$ is called
{\it right balanced} by $\mathscr{F}\times\mathscr{G}$ if for any $M\in\mathscr{C}$, there exists a
$T(-,\mathscr{G})$-exact complex $\cdots\to F_1\to F_0\to M\to 0$ with each $F_i\in\mathscr{F}$,
and for any $N\in\mathscr{D}$, there exists a $T(\mathscr{F},-)$-exact complex $0\to N\to G^0\to G^1\to\cdots$
with each $G^i\in\mathscr{G}$. They showed that if $T$ is right balanced by $\mathscr{F}\times\mathscr{G}$, and
if $F_\bullet\to M$ is a $T(-,\mathscr{G})$-exact complex and $N\to G^\bullet$ is a $T(\mathscr{F},-)$-exact
complex, then the complexes $T(F_\bullet,N)$ and $T(M,G^\bullet)$ have isomorphic homology. There are many
examples of right balanced functors in the module category when we regard $T$ as $\Hom$, see
[EJ2, Chapter 8]. Recently, Chen introduced in [C1] the notion of balanced pairs of additive subcategories
in an abelian category. Let $\mathscr{A}$ be an abelian category with
enough projectives and injectives. We use $\mathscr{P}(\mathscr{A})$ and $\mathscr{I}(\mathscr{A}$) to denote
the full subcategories of $\mathscr{A}$ consisting of projectives and injectives respectively. It is known
that the pair ($\mathscr{P}(\mathscr{A})$, $\mathscr{I}(\mathscr{A}$)) is a balanced pair, which is called
the {\it classical balanced pair}. Chen showed that for a balanced pair ($\mathscr{X}$, $\mathscr{Y}$)
of $\mathscr{A}$, it inherits some nice properties from the classical one.

The notion of cotorsion pairs was first introduced by Salce in [S], and it has been deeply studied in
homological algebra, representation theory and triangulated categories in recent years, see [EJTX], [Ho], [HI],
[IY], [KS], and so on. In particular, Hovey established in [Ho] a connection
between cotorsion pairs in abelian categories and model category theory.
In classical homological algebra, the definition of the cotorsion pair is based on the functor
$\Ext^i_\mathscr{A}(-,-)$. The advantage is that this functor is independent of the choices of the
projective resolutions of the first variable, and also independent of the choices of the injective resolutions
of the second variable. In other words, the cotorsion pair is essentially based on the balanced pair
($\mathscr{P}(\mathscr{A})$, $\mathscr{I}(\mathscr{A}$)). Based on these backgrounds mentioned above,
it is natural for us to introduce and study cotorsion pairs relative to balanced pairs, and we show
that relative cotorsion pairs share many nice properties of the classical one. This paper is organized
as follows.

In Section 2, we give some terminology and some preliminary results.

In Section 3, for an abelian category $\mathscr{A}$, we introduce the notion of cotorsion pairs relative
to a given balanced pair ($\mathscr{X}$, $\mathscr{Y}$). Similar to the classical case,
we also introduce the notions of complete, hereditary and perfect cotorsion pairs relative to
($\mathscr{X}$, $\mathscr{Y}$), and obtain some equivalent characterizations for the cotorsion pair
relative to ($\mathscr{X}$, $\mathscr{Y}$) being complete, hereditary and perfect respectively.

In Section 4, for a given balanced pair ($\mathscr{X}$, $\mathscr{Y}$) of the abelian category $\mathscr{A}$,
we introduce the notions of the right $\mathscr{X}$-derived category $D_{\mathbb{R}\mathscr{X}}^{\ast}$($\mathscr{A})$
and the left $\mathscr{Y}$-derived category $D_{\mathbb{L}\mathscr{Y}}^{\ast}$($\mathscr{A})$ of $\mathscr{A}$
for $\ast\in\{{\rm blank},-,+,b\}$. We show that in the bounded case, they are actually the same, and we denote both by
$D_\ast^b(\mathscr{A})$. Let ($\mathscr{X}$, $\mathscr{Y}$) be an admissible balanced pair. We give some criteria
for computing the $\mathscr{X}$-resolution dimension and the $\mathscr{Y}$-coresolution dimension
of an object in $\mathscr{A}$ in terms of the vanishing of relative cohomology groups. Moreover,
we show that if the $\mathscr{X}$-resolution dimension
of $\mathscr{Y}$ (resp. $\mathscr{Y}$-coresolution dimension
of $\mathscr{X}$) is finite, then the bounded homotopy category of $\mathscr{Y}$ (resp. $\mathscr{X}$)
is contained in that of $\mathscr{X}$ (resp. $\mathscr{Y}$). This generalizes a classical result of Happel.
As a consequence, we get that the right
$\mathscr{X}$-singularity category coincides with the left $\mathscr{Y}$-singularity category if
the $\mathscr{X}$-resolution dimension of $\mathscr{Y}$ and the $\mathscr{Y}$-coresolution dimension
of $\mathscr{X}$ are finite.

\vspace{0.5cm}

\centerline{\bf 2. Preliminaries}

\vspace{0.2cm}
Throughout this paper, $\mathscr{A}$ is an abelian category. For a subcategory of $\mathscr{A}$ we mean a full additive subcategory
closed under isomorphisms and direct summands. We use $\mathscr{P}(\mathscr{A})$ and $\mathscr{I}(\mathscr{A}$) to denote the
subcategories of $\mathscr{A}$ consisting of projective and injective objects respectively. We use $C(\mathscr{A})$ to denote the
category of complexes of objects in $\mathscr{A}$, $K^{*}(\mathscr{A})$ to denote the homotopy category of $\mathscr{A}$, and
$D^{*}(\mathscr{A})$ to denote the usual derived category by inverting the quasi-isomorphisms in $K^{*}(\mathscr{A})$, where
$* \in \{{\rm blank},-,+,b\}$.

Let $$X^{\bullet}:=\cdots \longrightarrow X^{n-1}\buildrel {d^{n-1}_{X}} \over \longrightarrow X^{n}\buildrel {d^{n}_{X}}
\over\longrightarrow X^{n+1} \to \cdots$$ be a complex in $C(\mathscr{A})$ and $f: X^{\bullet}\to Y^{\bullet}$ a cochain map in $C(\mathscr{A})$.
We use $\Con(f)$ to denote the mapping cone of $f$. Recall that $X^{\bullet}$ is called {\it acyclic} (or {\it exact}) if $H^{i}(X^{\bullet})=0$
for any $i\in \mathbb{Z}$ (the ring of integers), and $f$ is called a {\it quasi-isomorphism} if $H^{i}(f)$ is an isomorphism for any $i\in \mathbb{Z}$.
We have that $f$ is a quasi-isomorphism if and only if $\Con(f)$ is acyclic.

\vspace{0.2cm}

{\bf Definition 2.1.} (1) ([AR]) Let $\mathscr{X}\subseteq\mathscr{Y}$ be
subcategories of $\mathscr{A}$. A morphism $f: X\to Y$ in
$\mathscr{A}$ with $X\in\mathscr{X}$ and $Y\in\mathscr{Y}$ is called a {\it right $\mathscr{X}$-approximation} of $Y$
if for any morphism $g: X' \to Y$ in $\mathscr{A}$ with $X'\in\mathscr{X}$, there exists a morphism $h: X'\to X$
such that the following diagram commutes:
$$\xymatrix{ & X' \ar[d]^{g} \ar@{-->}[ld]_{h}\\
X \ar[r]^{f} & Y.}$$
If any endomorphism $s:X\to X$ is an automorphism whenever $f=fs$, then $f$ is called {\it right minimal}.
If each object in $\mathscr{Y}$ has a right $\mathscr{X}$-approximation,
then $\mathscr{X}$ is called {\it contravariantly finite} in $\mathscr{Y}$. Dually, the notions of
{\it left $\mathscr{X}$-approximations}, {\it left minimal morphisms} and {\it covariantly finite subcategories} are defined.

(2) ([C1]) A contravariantly finite subcategory $\mathscr{X}$ of $\mathscr{A}$ is called {\it admissible} if each right
$\mathscr{X}$-approximation is epic. Dually, the notion of {\it coadmissible subcategories} is defined.

\vspace{0.2cm}

{\bf Definition 2.2.}
(1) ([C1]) Given a subcategory $\mathscr{X}$ of $\mathscr{A}$. A complex $A^\bullet$ in $C(\mathscr{A})$ is called {\it right} (resp. {\it left})
{\it $\mathscr{X}$-acyclic} if the complex $\Hom_\mathscr{A}(X,A^\bullet)$ (resp. $\Hom_\mathscr{A}(A^\bullet,X)$ ) is acyclic for any
$X\in\mathscr{X}$. A cochain map $f:A^\bullet\to B^\bullet$ in $C(\mathscr{A})$ is said to be {\it right} (resp. {\it left)
$\mathscr{X}$-quasi-isomorphism} if the cochain (resp. chain) map $\Hom_\mathscr{A}(X,f)$ (resp. $\Hom_\mathscr{A}(f,X)$) is a quasi-isomorphism
for any $X\in\mathscr{X}$. It is equivalent to that $\Con(f)$ is right (resp. left) $\mathscr{X}$-acyclic.

(2) ([C1, EJ2]) Given a contravariantly finite subcategory $\mathscr{X}$ of $\mathscr{A}$ and an object $M\in \mathscr{A}$.
An {\it $\mathscr{X}$-resolution} of $M$ is a complex $$\cdots\to X^{-2}\buildrel {d^{-2}} \over \longrightarrow X^{-1}\buildrel
{d^{-1}} \over \longrightarrow X^0 \buildrel {\varepsilon} \over \longrightarrow M\to 0$$ in $\mathscr{A}$ with each $X^i\in \mathscr{X}$
such that it is right $\mathscr{X}$-acyclic. Usually we denote the complex by $X^\bullet\buildrel {\varepsilon} \over \longrightarrow M $
for short, where $$X^\bullet:=\cdots\to X^{-2}\buildrel {d^{-2}} \over \longrightarrow X^{-1}\to \cdots\buildrel {d^{-1}} \over \longrightarrow X^0\to 0$$
is the deleted $\mathscr{X}$-resolution of $M$. The {\it $\mathscr{X}$-resolution dimension} $\mathscr{X}$-$\resdim M$ of $M$
is defined to be the minimal integer $n\geq 0$ such that there exists an $\mathscr{X}$-resolution:
$$0\to X^{-n}\to \cdots\to X^{-1}\to X^0\to M\to 0.$$ If no such an integer exists, we set $\mathscr{X}$-$\resdim M =\infty$.
The {\it global $\mathscr{X}$-resolution dimension} $\mathscr{X}$-$\resdim \mathscr{A}$ of $\mathscr{A}$ is defined to be the
supreme of the $\mathscr{X}$-resolution dimensions of all objects in $\mathscr{A}$.

Dually, if $\mathscr{X}$ is a covariantly finite subcategory of $\mathscr{A}$, then the notions of {\it $\mathscr{X}$-coresolutions},
{\it $\mathscr{X}$-coresolution dimensions} and the {\it global $\mathscr{X}$-coresolution
dimension} are defined.

\vspace{0.2cm}

{\bf Lemma 2.3.}  {\it Let $\mathscr{X}$ be a subcategory of $\mathscr{A}$ and let

$$0\to L\buildrel {f} \over \longrightarrow M\buildrel {g} \over \longrightarrow N\to 0\eqno{(2.1)}$$ be an acyclic complex.

(1) If (2.1) is right $\mathscr{X}$-acyclic,
then for any morphisms $N'\buildrel {\alpha} \over \longrightarrow N$ and $L\buildrel {s} \over \longrightarrow L'$,
we have the following pull-back diagram with the upper row right $\mathscr{X}$-acyclic:
$$\xymatrix{0\ar[r]&L\ar[r]^{f'}\ar@{=}[d]&M'\ar[r]^{g'}\ar[d]^{\beta}&N'\ar[d]^\alpha\ar[r]&0\\
0\ar[r]&L\ar[r]^{f}&M\ar[r]^{g}&N\ar[r]&0,
}$$ and the following push-out diagram with the bottom row right $\mathscr{X}$-acyclic:
$$\xymatrix{0\ar[r]&L\ar[r]^{f}\ar[d]^s&M\ar[r]^{g}\ar[d]^{t}&N\ar@{=}[d]\ar[r]&0\\
0\ar[r]&L'\ar[r]^{f''}&M''\ar[r]^{g''}&N\ar[r]&0.
}$$

(2) If (2.1) is left $\mathscr{X}$-acyclic,
then for any morphisms $N'\buildrel {\alpha} \over \longrightarrow N$ and $L\buildrel {s} \over \longrightarrow L'$,
we have the following pull-back diagram with the upper row left $\mathscr{X}$-acyclic:
$$\xymatrix{0\ar[r]&L\ar[r]^{f'}\ar@{=}[d]&M'\ar[r]^{g'}\ar[d]^{\beta}&N'\ar[d]^\alpha\ar[r]&0\\
0\ar[r]&L\ar[r]^{f}&M\ar[r]^{g}&N\ar[r]&0,
}$$ and the following push-out diagram with the bottom row left $\mathscr{X}$-acyclic:
$$\xymatrix{0\ar[r]&L\ar[r]^{f}\ar[d]^s&M\ar[r]^{g}\ar[d]^{t}&N\ar@{=}[d]\ar[r]&0\\
0\ar[r]&L'\ar[r]^{f''}&M''\ar[r]^{g''}&N\ar[r]&0.
}$$}

\vspace{0.2cm}

{\it Proof.} (1) Because the sequence (2.1) is right $\mathscr{X}$-acyclic by assumption, for any morphism $h:X\to N'$
with $X\in \mathscr{X}$ there exists a morphism $i:X\to M$ such that $\alpha h=gi$. Since the right square in the first diagram
is a pull-back diagram, there exists a morphism $\phi:X\to M'$  such that $h=g'\phi$. It implies that the upper row in this
diagram is right $\mathscr{X}$-acyclic.

Also because the sequence (2.1) is right $\mathscr{X}$-acyclic acyclic, for any morphism $h':X\to N$ with $X\in \mathscr{X}$
there exists a morphism $i':X\to M$ such that $h'=gi'=g''ti'$. It implies that the bottom row in the second
diagram is right $\mathscr{X}$-acyclic.

(2) It is dual to (1). \hfill{$\square$}

\vspace{0.2cm}

{\bf Lemma 2.4.} {\it Let $A^{\bullet}$ be a complex in $C(\mathscr{A})$. Then $A^{\bullet}$ is right $\mathscr{X}$-acyclic
if and only if the complex $\Hom_\mathscr{A}(X^{\bullet},A^{\bullet})$ is acyclic for any $X^{\bullet}\in K^{-}(\mathscr{X})$.}

\vspace{0.2cm}

{\it Proof.} See [CFH, Lemma 2.4]. \hfill{$\square$}

\vspace{0.2cm}

{\bf Lemma 2.5.} {\it (1) Let $X^{\bullet}$ be a complex in $K^{-}(\mathscr{X})$ and let $f$ : $A^{\bullet}\to X^{\bullet}$ be
a right $\mathscr{X}$-quasi-isomorphism in $C(\mathscr{A})$. Then there exists a cochain map $g: X^{\bullet}\to A^{\bullet}$
such that $fg$ is homotopic to $\id_{X^{\bullet}}$.

(2) Any right $\mathscr{X}$-quasi-isomorphism between two complexes in $K^{-}(\mathscr{X})$ is a homotopy equivalence.}

\vspace{0.2cm}
{\it Proof.} (1) Consider the distinguished triangle:
$$A^{\bullet}\buildrel {f} \over \longrightarrow X^{\bullet}\to
\Con(f)\to A^{\bullet}[1]$$ in $K(\mathscr{A})$ with $\Con(f)$ right $\mathscr{X}$-acyclic. By applying the functor
$\Hom_{K(\mathscr{A})}(X^\bullet,-)$ to it, we get an exact sequence:
$$\Hom_{K(\mathscr{A})}(X^\bullet,A^\bullet)\buildrel {\Hom_{K(\mathscr{A})}(X^\bullet,f)} \over \longrightarrow
\Hom_{K(\mathscr{A})}(X^\bullet,X^\bullet)\to \Hom_{K(\mathscr{A})}(X^\bullet,\Con(f)).$$
It follows from Lemma 2.4 that $\Hom_{K(\mathscr{A})}(X^\bullet,\Con(f))\cong H^0\Hom_\mathscr{A}(X^{\bullet},\Con(f))=0$.
So there exists a cochain map $g: X^{\bullet}\to A^{\bullet}$ such that $fg$ is homotopic to $\id_{X^{\bullet}}$.

(2) It is a consequence of (1).  \hfill{$\square$}

\vspace{0.5cm}

\centerline{\bf 3. Cotorsion pairs relative to balanced pairs}

\vspace{0.2cm}

{\bf Definition 3.1.} (see [C1] and [EJ2]) A pair ($\mathscr{X}$, $\mathscr{Y}$) of subcategories of $\mathscr{A}$ is called a {\it balanced pair} if the following
conditions are satisfied:

(1) $\mathscr{X}$ is contravariantly finite in $\mathscr{A}$ and $\mathscr{Y}$ is covariantly finite in $\mathscr{A}$.

(2) For any object $M\in\mathscr{A}$, there exists an $\mathscr{X}$-resolution $X^\bullet\longrightarrow M $ of $M$ such that it is left $\mathscr{Y}$-acyclic.

(3) For any object $N\in \mathscr{A}$, there exists a $\mathscr{Y}$-coresolution $N\longrightarrow Y^\bullet $ of $N$ such that it is right $\mathscr{X}$-acyclic.

\vspace{0.2cm}

We list some examples of balanced pairs as follows.

\vspace{0.2cm}

{\bf Example 3.2.} (1) Recall that $\mathscr{A}$ is said to have {\it enough projectives} (resp. {\it enough injectives}) if for any $M\in \mathscr{A}$,
there exists an epimorphism $P\to M\to 0$ (resp. a monomorphism $0\to M\to I$) with $P$ (resp. $I$) in $\mathscr{P}(\mathscr{A})$ (resp. $\mathscr{I}(\mathscr{A})$).
In case for $\mathscr{A}$ having enough projectives and injectives, it is well known that the pair ($\mathscr{P}(\mathscr{A})$, $\mathscr{I}(\mathscr{A}$))
is a balanced pair. We call it the {\it classical balanced pair}.

(2) ([EJ2, Example 8.3.2]) Let $R$ be a ring and $\Mod R$ the category of left $R$-modules, and let $\mathscr{PP}(R)$ and
$\mathscr{PI}(R)$ be the subcategories of $\Mod R$ consisting of pure projective modules and pure injective modules respectively.
Then ($\mathscr{PP}(R),\mathscr{PI}(R)$) is a balanced pair in $\Mod R$.

(3) ([EJ2, Theorem 12.1.4]) Let $R$ be an $n$-Gorenstein ring (that is, $R$ is a left and right Noetherian ring with left and right self-injective dimensions at most $n$),
and let $\GProj R$ and $\GInj R$ be the subcategories of $\Mod R$ consisting of Gorenstein projective and Gorenstein injective modules respectively.
Then the pair ($\GProj R$, $\GInj R$) is a balanced pair in $\Mod R$.

\vspace{0.2cm}

Let $\mathscr{X}$ (resp. $\mathscr{Y}$) be a contravariantly finite (resp. covariantly finite) subcategory of $\mathscr{A}$. Then the pair
($\mathscr{X}$, $\mathscr{Y}$) is a balanced pair if and only if the class of right $\mathscr{X}$-acyclic complexes
coincides with that of left $\mathscr{Y}$-acyclic complexes ([C1, Proposition 2.2]). In what follows, we call a complex {\it $\ast$-acyclic}
if it is both right $\mathscr{X}$-acyclic and left $\mathscr{Y}$-acyclic.

\vspace{0.2cm}

{\bf Definition 3.3.} ([EJ2, Definition 7.1.2]) Let $\mathscr{A}$ have enough projectives and injectives.
A pair ($\mathscr{C}$, $\mathscr{D}$) of subcategories of $\mathscr{A}$ is called a {\it cotorsion pair}
if $\mathscr{C}={^\perp\mathscr{D}}$ and $\mathscr{D}=\mathscr{C}^\perp$, where $^\perp\mathscr{D}=\{C\in\mathscr{A}\mid\Ext_\mathscr{A}^1(C,D)=0$ for any $D\in\mathscr{D}$\}
and $\mathscr{C}^\perp=\{D\in\mathscr{A}\mid\Ext_\mathscr{A}^1(C,D)=0$ for any $C\in\mathscr{C}$\}.

\vspace{0.2cm}

Notice that the functor $\Ext_\mathscr{A}^1(-,-)$ is based on the classical balanced pair ($\mathscr{P}(\mathscr{A})$, $\mathscr{I}(\mathscr{A})$),
it induces an isomorphism of cohomology groups whether we take a projective resolution of the first variable or take an injective coresolution of the second variable.
From this viewpoint we may say that the cotorsion pair defined above is {\it a cotorsion pair} relative to the balanced pair ($\mathscr{P}(\mathscr{A})$, $\mathscr{I}(\mathscr{A})$).

Let ($\mathscr{X}$, $\mathscr{Y}$) be a balanced pair and $M,N\in \mathscr{A}$. Choose an $\mathscr{X}$-resolution $X^\bullet\longrightarrow M$ of $M$ and
a $\mathscr{Y}$-coresolution $N\longrightarrow Y^\bullet$ of $N$. We get two cohomological groups $\Ext_\mathscr{X}^i(M,N):=H^i(\Hom_{\mathscr{A}}(X^\bullet,N))$ and
$\Ext_\mathscr{Y}^i(M,N):=H^i(\Hom_{\mathscr{A}}(M,Y^\bullet))$ for any $i\in\mathbb{Z}$. They are independent of the choices of the $\mathscr{X}$-resolutions of $M$ and the $\mathscr{Y}$-coresolutions of $N$
respectively. For any $i\in \mathbb{Z}$, there exists an isomorphism of abelian groups $\Ext_\mathscr{X}^i(M,N)\cong \Ext_\mathscr{Y}^i(M,N)$ ([EJ2]). We denote both abelian groups by $\Ext_\ast^i(M,N)$.
Motivated by the above argument, we introduce the following

\vspace{0.2cm}

{\bf Definition 3.4.} Let ($\mathscr{X}$, $\mathscr{Y}$) be a balanced pair. A pair ($\mathscr{C}$, $\mathscr{D}$) of subcategories of $\mathscr{A}$ is called a {\it cotorsion pair}
relative to ($\mathscr{X}$, $\mathscr{Y}$) if $\mathscr{C}$=$^{\perp_\ast}\mathscr{D}$ and $\mathscr{D}$=$\mathscr{C}^{\perp_\ast}$, where
$^{\perp_\ast}\mathscr{D}=\{C\in\mathscr{A}\mid\Ext_\ast^1(C,D)=0$ for any $D\in\mathscr{D}\}$ and $\mathscr{C}^{\perp_\ast}=\{D\in\mathscr{A}\mid\Ext_\ast^1(C,D)=0$
for any $C\in\mathscr{C}\}$.

\vspace{0.2cm}

{\it In the rest of this section, we fix a balanced pair ($\mathscr{X}$, $\mathscr{Y}$) and a cotorsion pair ($\mathscr{C}$, $\mathscr{D}$) relative to ($\mathscr{X}$, $\mathscr{Y}$) in $\mathscr{A}$.}

\vspace{0.2cm}

{\bf Proposition 3.5.} {\it For any $A\in \mathscr{A}$, we have $\Ext_\ast^{\geq 1}(X,A)=0=\Ext_\ast^{\geq 1}(A,Y)$ for any $X\in\mathscr{X}$ and $Y\in\mathscr{Y}$.}

\vspace{0.2cm}

{\it Proof.} It is straightforward. \hfill{$\square$}

\vspace{0.2cm}

{\bf Definition 3.6.} (1) Let $\mathscr{E}$ be a subcategory of $\mathscr{A}$. $\mathscr{E}$ is said to be {\it closed under $\ast$-extensions}
if for any $\ast$-acyclic complex $0\to L\to M\to N\to 0$ in $\mathscr{A}$, $L,N\in \mathscr{E}$ implies $M\in \mathscr{E}$; $\mathscr{E}$ is said to be {\it closed under $\ast$-epimorphisms}
if for any $\ast$-acyclic complex $0\to L\to M\to N\to 0$ in $\mathscr{A}$, $M,N\in \mathscr{E}$ implies $L\in \mathscr{E}$; $\mathscr{E}$ is said to be {\it closed under $\ast$-monomorphisms}
if for any $\ast$-acyclic complex $0\to L\to M\to N\to 0$ in $\mathscr{A}$, $L,M\in \mathscr{E}$ implies $N\in \mathscr{E}$.

(2) A subcategory $\mathscr{E}$ of $\mathscr{A}$ is called {\it $\mathscr{X}$-resolving} if $\mathscr{X}\subseteq\mathscr{E}$ and $\mathscr{E}$ is closed under
$\ast$-extensions and $\ast$-epimorphisms; and $\mathscr{E}$ is called {\it $\mathscr{Y}$-coresolving} if $\mathscr{Y}\subseteq\mathscr{E}$ and $\mathscr{E}$ is closed under $\ast$-extensions
and $\ast$-monomorphisms.

\vspace{0.2cm}

{\bf Proposition 3.7.}

{\it (1) $\mathscr{X}\subseteq\mathscr{C}$ and $\mathscr{Y}\subseteq\mathscr{D}$.

(2) Both $\mathscr{C}$ and $\mathscr{D}$ are closed under $\ast$-extensions.}

\vspace{0.2cm}

{\it Proof.} (1) It is trivial.

(2) Let $$0\to L\to M\to N\to 0$$
be a $\ast$-acyclic complex in $\mathscr{A}$ with $L,N\in\mathscr{C}$. By [EJ2, Theorem 8.2.3], for any $D\in\mathscr{D}$ we have an acyclic complex
$$\Ext_\ast^1(N,D)\to\Ext_\ast^1(M,D)\to\Ext_\ast^1(L,D).$$ Then we have $\Ext_\ast^1(N,D)=0=\Ext_\ast^1(L,D)$.
So $\Ext_\ast^1(M,D)=0$ and $M\in\mathscr{C}$. Thus $\mathscr{C}$ is closed under $\ast$-extensions. Similarly, we have that $\mathscr{D}$ is closed under $\ast$-extensions.
\hfill{$\square$}

\vspace{0.2cm}

The following result is a relative version of [AR, Lemmas 3.1 and 3.2].

\vspace{0.2cm}

{\bf Theorem 3.8.}  {\it The following statements are equivalent.

(1) $\mathscr{C}$ is $\mathscr{X}$-resolving.

(2) $\mathscr{D}$ is $\mathscr{Y}$-coresolving.

(3) $\Ext_\ast^{\geq 1}(C,D)=0$ for any $C\in\mathscr{C}$ and $D\in\mathscr{D}$.

In this case, ($\mathscr{C}$, $\mathscr{D}$) is called hereditary.}

\vspace{0.2cm}

{\it Proof.} $(1)\Rightarrow (3)$. Let $C\in\mathscr{C}$ and $D\in\mathscr{D}$, and let $$0\to K\to X\to C\to 0$$ be a $\ast$-acyclic complex in $\mathscr{A}$
with $X\in\mathscr{X}$. Then $X\in\mathscr{C}$ by Proposition 3.7(1), and so $K\in\mathscr{C}$ by (1). Hence $\Ext_\ast^1(K,D)=0$. By [EJ2, Theorem 8.2.3], we have an exact sequence:
$$\Ext_\ast^1(X,D)\to \Ext_\ast^1(K,D)\to\Ext_\ast^2(C,D)\to\Ext_\ast^2(X,D).$$
Since $\Ext_\ast^1(X,D)=0=\Ext_\ast^2(X,D)$ by Proposition 3.5, we have $\Ext_\ast^2(C,D)\cong \Ext_\ast^1(K,D)=0$. We get $\Ext_\ast^{i\geq 1}(C,D)=0$ inductively.

$(3)\Rightarrow (1)$. Let $$0\to L\to M\to N\to 0$$ be a $\ast$-acyclic complex in $\mathscr{A}$ with $M,N\in \mathscr{C}$. By [EJ2, Theorem 8.2.3], for any $D\in \mathscr{D}$
we have the following exact sequence:
$$\Ext_\ast^1(M,D)\to \Ext_\ast^1(L,D)\to\Ext_\ast^2(N,D).$$
Because $\Ext_\ast^1(M,D)=0=\Ext_\ast^2(N,D)$ by assumption, we have $\Ext_\ast^1(L,D)=0$ and $L\in \mathscr{C}$. Now the assertion follows from Proposition 3.7.

Dually, we get $(2)\Leftrightarrow (3)$. \hfill{$\square$}

\vspace{0.2cm}

{\bf Definition 3.9.} (cf. [EJ2, Definition 7.1.5])

(1) ($\mathscr{C}$, $\mathscr{D}$) is said to have {\it enough projectives} if for any $M\in\mathscr{A}$,
there exists a $\ast$-acyclic complex
$$0\to D\to C\to M\to 0$$ in $\mathscr{A}$ with $C\in\mathscr{C}$ and $D\in\mathscr{D}$; and it is said to have {\it enough injectives}
if for any $M\in\mathscr{A}$, there exists a
$\ast$-acyclic complex $$0\to M\to D\to C\to 0$$ in $\mathscr{A}$ with $C\in\mathscr{C}$ and $D\in\mathscr{D}$.

(2) If ($\mathscr{C}$, $\mathscr{D}$) has enough projectives and enough injectives, then it is called {\it complete}.

\vspace{0.2cm}

{\bf Proposition 3.10.} {\it $\mathscr{X}$ is admissible if and only if $\mathscr{Y}$ is coadmissible.
In this case, ($\mathscr{X}$, $\mathscr{Y}$) is called admissible.}

\vspace{0.2cm}

{\it Proof.} See [C1, Corollary 2.3]. \hfill{$\square$}

\vspace{0.2cm}

It is obvious that if ($\mathscr{X}$, $\mathscr{Y}$) is admissible, then each $\ast$-acyclic complex is acyclic.

\vspace{0.2cm}

{\bf Theorem 3.11.} {\it If ($\mathscr{X}$, $\mathscr{Y}$) is admissible, then ($\mathscr{C}$, $\mathscr{D}$) has enough projectives if and only if it has enough injectives.}
\vspace{0.2cm}

{\it Proof.} We only show the ``if" part, and the ``only if" part follows dually.

Assume that ($\mathscr{C}$, $\mathscr{D}$) has enough injectives and $M\in\mathscr{A}$. Choose a $\ast$-acyclic complex
$$0\to K\to X\to M\to 0$$ in $\mathscr{A}$ with $X\in\mathscr{X}$. Since ($\mathscr{C}$, $\mathscr{D}$) has enough injectives, there exists a $\ast$-acyclic complex
$$0\to K\to D\to C\to 0$$ in $\mathscr{A}$ with $C\in\mathscr{C}$ and $D\in\mathscr{D}$. Because ($\mathscr{X}$, $\mathscr{Y}$) is admissible, each $\ast$-acyclic complex is acyclic.
So we have the following push-out diagram with acyclic columns and rows:
$$\xymatrix{ &0\ar[d]&0\ar[d]& & \\
0\ar[r]&K\ar[r]\ar[d]&X\ar[r]\ar[d]&M\ar[r]\ar@{=}[d]&0\\
0\ar[r]&D\ar[r]\ar[d]&E\ar[r]\ar[d]&M\ar[r]&0\\
 &C\ar[d]\ar@{=}[r]&C\ar[d]& & \\
 &0&0.& &
}$$
Then all of columns and rows are $\ast$-acyclic by Lemma 2.3. Since $X,C\in\mathscr{C}$, it follows from Proposition 3.7 that $E\in\mathscr{C}$.
The assertion follows. \hfill{$\square$}

\vspace{0.2cm}

{\bf Lemma 3.12.} {\it Let ($\mathscr{X}$, $\mathscr{Y}$) be admissible and $\mathscr{E}$ a subcategory of $\mathscr{A}$ which is closed under $\ast$-extensions.

(1) If $\varphi: E\to M$ is a minimal right $\mathscr{E}$-approximation and
$$\xymatrix{0\ar[r]&S\ar[d]^f\ar[r]^i&P\ar[d]^\theta\ar[r]^\pi&G\ar[r]&0\\
&E\ar[r]^\varphi&M & &
}$$
is a commutative diagram with $G\in\mathscr{E}$ such that the upper row is $\ast$-acyclic, then there exists a morphism $\alpha:P\to E$ such that $f=\alpha i$ and $\theta=\varphi\alpha$.

(2) If $\psi: M\to E$ is a minimal left $\mathscr{E}$-approximation and
$$\xymatrix{&&M\ar[r]^\psi\ar[d]^\theta&E\ar[d]^f&\\
0\ar[r]&F\ar[r]^i&Q\ar[r]^\pi&K\ar[r]&0
}$$
is a commutative diagram with $F\in\mathscr{E}$ such that the bottom row is $\ast$-acyclic, then there exists a morphism $\alpha:E\to Q$ such that $\theta=\alpha \psi$ and $f=\pi\alpha$.}

\vspace{0.2cm}

{\it Proof.} (1) Consider the following push-out diagram:
$$\xymatrix{0\ar[r]&S\ar[d]^f\ar[r]^i&P\ar[d]^k\ar[r]^\pi&G\ar[r]\ar@{=}[d]&0\\
0\ar[r]&E\ar[r]^j&X\ar[r]&G\ar[r]&0.
}$$
Because the first row is $\ast$-acyclic by assumption, it follows from Lemma 2.3 that the bottom row is also $\ast$-acyclic.
Since $E,G\in\mathscr{E}$, we have $X\in\mathscr{E}$. By the universal property of push-outs there exists a morphism $h:X\to M$ such that $\varphi=hj$ and $\theta=hk$.
Because $\varphi:E\to M$ is a minimal right $\mathscr{E}$-approximation by assumption, there exists a morphism $g:X\to E$ such that $h=\varphi g$.
Thus $\varphi=hj=\varphi gj$, which implies that $gj:E\to E$ is an automorphism. We may assume $gj=\id_E$. Then by letting $\alpha=gk$,
we have $jf=ki=jgki=j\alpha i$. Since $j$ is a monomorphism, $f=\alpha i$. It follows from $\theta=hk$ and $h=\varphi g$ that $\theta=hk=\varphi gk=\varphi\alpha$, we complete the proof.

(2) It is dual to (1). \hfill{$\square$}

\vspace{0.2cm}

The following result is a relative version of the Wakamatsu's lemma.

\vspace{0.2cm}

{\bf Proposition 3.13.} {\it Let ($\mathscr{X}$, $\mathscr{Y}$) be admissible and $\mathscr{E}$ a subcategory of $\mathscr{A}$ which is closed under $\ast$-extensions.
Then we have

(1) The kernel of every minimal right $\mathscr{E}$-approximation is in $\mathscr{E}^{\perp_\ast}$.

(2) The cokernel of every minimal left $\mathscr{E}$-approximation is in $^{\perp_\ast}\mathscr{E}$.
 }

\vspace{0.2cm}

{\it Proof.} (1) Let $\varphi: E\to M$ be a minimal right $\mathscr{E}$-approximation of an object $M$ in $\mathscr{A}$, and let $K:=\Ker\varphi$ and $i:K\to E$ be the inclusion.
Because $\mathscr{X}$ is contravariantly finite in $\mathscr{A}$, for any $E'\in \mathscr{E}$ there exists a $\ast$-acyclic complex
$$0\to S\to X\to E'\to 0$$ in $\mathscr{A}$ with $X\in\mathscr{X}$. By applying the functor $\Hom_\mathscr{A}(-,K)$ we get an exact sequence:
$$\Hom_\mathscr{A}(X,K)\to \Hom_\mathscr{A}(S,K)\to\Ext^1_\ast(E',K)\to 0.$$ For any morphism $f:S\to K$, it follows from Lemma 3.12 that there exists a morphism
$g:X\to E$ such that the following diagram
$$\xymatrix{S\ar[r]\ar[d]_{if}&X\ar[d]^0\ar@{-->}[ld]_{g}\\
E\ar[r]^\varphi&M
}$$ is commutative.
Then $\Im g\subseteq \Ker\varphi=K$. So the map $\Hom_\mathscr{A}(X,K)\to \Hom_\mathscr{A}(S,K)$ is epic and $\Ext^1_\ast(E',K)=0$.

(2) It is dual to (1). \hfill{$\square$}

\vspace{0.2cm}

{\bf Definition 3.14.} ($\mathscr{C}$, $\mathscr{D}$) is called {\it perfect} if every object of $\mathscr{A}$ has a minimal right
$\mathscr{C}$-approximation and a minimal left $\mathscr{D}$-approximation.

\vspace{0.2cm}

Let ($\mathscr{X}$, $\mathscr{Y}$) be admissible. If ($\mathscr{C}$, $\mathscr{D}$) is perfect, then it is complete by Proposition 3.13.
The following result is a relative version of [EJTX, Theorem 3.8].

\vspace{0.2cm}

{\bf Theorem 3.15.}  {\it Let ($\mathscr{X}$, $\mathscr{Y}$) be admissible and ($\mathscr{C}$, $\mathscr{D}$) a hereditary cotorsion pair.
Then the following statements are equivalent.

(1) ($\mathscr{C}$, $\mathscr{D}$) is perfect.

(2) Every object of $\mathscr{A}$ has a minimal right $\mathscr{C}$-approximation and every object of $\mathscr{C}$ has a minimal left $\mathscr{D}$-approximation.

(3) Every object of $\mathscr{A}$ has a minimal left $\mathscr{D}$-approximation and every object of $\mathscr{D}$ has a minimal right $\mathscr{C}$-approximation.}

\vspace{0.2cm}

{\it Proof.} Both $(1)\Rightarrow (2)$ and $(1)\Rightarrow (3)$ are trivial. In the following we only prove $(2)\Rightarrow (1)$, and $(3)\Rightarrow (1)$ follows dually.

Let $\varphi:C\to M$ be a minimal right $\mathscr{C}$-approximation of an object $M$ in $\mathscr{A}$. Since ($\mathscr{X}$, $\mathscr{Y}$) is admissible, by Proposition 3.13(1)
that there exists a $\ast$-acyclic complex: $$0\to D\buildrel {i} \over \longrightarrow C\buildrel {\varphi} \over \longrightarrow M\to 0$$ in $\mathscr{A}$ with
$D\in\mathscr{D}$. Let $\psi:C\to D'$ be a minimal left $\mathscr{D}$-approximation of $C$. Then by Proposition 3.13(2), we get a $\ast$-acyclic complex
$$0\to C\buildrel {\psi} \over \longrightarrow D'\buildrel {\pi} \over \longrightarrow C'\to 0$$ in $\mathscr{A}$ with $C'\in\mathscr{C}$. Since $C,C'\in\mathscr{C}$,
we have $D'\in\mathscr{C}\cap\mathscr{D}$. Consider the following push-out diagram:
$$\xymatrix{ & & 0\ar[d]& 0\ar[d]& \\
0\ar[r] & D\ar@{=}[d]\ar[r]^i & C\ar[d]^\psi\ar[r]^\varphi& M \ar[d]^{\psi'}\ar[r]& 0\\
0\ar[r] & D\ar[r]^{i'} & D^{'}\ar[d]^\pi\ar[r]^{\varphi'}&X \ar[d]\ar[r]& 0\\
 & & C'\ar[d]\ar@{=}[r] & C'\ar[d]& \\
 & & 0 & 0.& }$$
By Lemma 2.3 the rightmost column is $\ast$-acyclic. For any morphism $f:D\to Y$ with $Y\in\mathscr{Y}$,
there exists a morphism $g:C\to Y$ such that $f=gi$, and then there exists a morphism $j:D'\to Y$ such that $g=j\psi$.
It follows that $f=gi=j\psi i=ji'$ and the middle row is $\ast$-acyclic. Because ($\mathscr{C}$, $\mathscr{D}$) is hereditary
by assumption, we have $X\in\mathscr{D}$.

To get the desired assertion, it suffices to show $\psi':M\to X$ is left minimal. Let $h:X\to X$ satisfying $\psi'=h\psi'$. By applying the functor $\Hom_\mathscr{A}(D',-)$
to the middle row, we have a morphism $h':D'\to D'$ such that the following diagram
$$\xymatrix{D'\ar[r]^{\varphi'}\ar@{-->}[d]^{h'} &X\ar[d]^{h} \\
D'\ar[r]^{\varphi'} &X }$$ is commutative.
Hence we have the following commutative diagram:
$$\xymatrix{D'\ar[r]^{\varphi'}\ar@{-->}[d]^{h'} &X\ar[d]^{h}&M\ar[l]_{\psi'}\ar@{=}[d] \\
D'\ar[r]^{\varphi'} &X&M\ar[l]_{\psi'}. }$$
Note that the diagram $$\xymatrix{C\ar[r]^{\varphi}\ar[d]^{\psi} &M\ar[d]^{\psi'} \\
D'\ar[r]^{\varphi'} &X }$$ is both a push-out diagram and a pull-back diagram. Then there exists morphism $h'':C\to C$ such that the following diagram
$$\xymatrix@!0{ & C\ar[rr]\ar@{-->}[dd]\ar[ld]&&M\ar[ld]\ar@{=}[dd]\\
D'\ar[rr]\ar[dd]&&X\ar[dd]\\
&C\ar[rr]\ar[ld]&&M\ar[ld]\\
D'\ar[rr]&&X
}$$ is commutative.
Since $\varphi:C\to M$ is right minimal, $h'':C\to C$ is an automorphism. Then it follows from the left minimality of $\psi:C\to D'$ that $h':D'\to D'$
is also an automorphism. It implies that $h:X\to X$ is an automorphism and $\psi':M\to X$ is left minimal.   \hfill{$\square$}

\vspace{0.5cm}

\centerline{\bf 4. Derived categories relative to balanced pairs}

\vspace{0.2cm}

Let $\mathscr{X}$ be a subcategory of $\mathscr{A}$. It is known that $K^{\ast}(\mathscr{A})$ is a triangulated category for $\ast \in \{{\rm blank},-,+,b\}$.
Denote by $K^{\ast}_{\mathbb{R}\mathscr{X}-ac}(\mathscr{A})$ (resp. $K^{\ast}_{\mathbb{L}\mathscr{X}-ac}(\mathscr{A})$) the full triangulated subcategory of
$K^{\ast}(\mathscr{A})$ consisting of right $\mathscr{X}$-acyclic (resp. left $\mathscr{X}$-acyclic) complexes. Both of them are thick subcategories because
they are closed under direct summands. Denote by $\sum_{\mathbb{R}\mathscr{X}}^\ast$ (resp. $\sum_{\mathbb{L}\mathscr{X}}^\ast$) the class of all right
(resp. left) $\mathscr{X}$-quasi-isomorphisms in $K^{\ast}(\mathscr{A})$. Then a cochain map is a right (resp. left) $\mathscr{X}$-quasi-isomorphism
if and only if its mapping cone is right (resp. left) $\mathscr{X}$-acyclic. Thus $\sum^\ast_{\mathbb{R}\mathscr{X}}$ (resp. $\sum_{\mathbb{L}\mathscr{X}}^\ast$)
is the saturated compatible multiplicative system determined by $K^{\ast}_{\mathbb{R}\mathscr{X}-ac}(\mathscr{A})$ (resp. $K^{\ast}_{\mathbb{L}\mathscr{X}-ac}(\mathscr{A})$).
\vspace{0.2cm}

{\bf Definition 4.1.} ([V]) The Verdier quotient category $D_{\mathbb{R}\mathscr{X}}^{\ast}$($\mathscr{A}):
=K^{\ast}(\mathscr{A})/K_{\mathbb{R}\mathscr{X}{\text -ac}}^{\ast}(\mathscr{A})$
is called the {\it right $\mathscr{X}$-derived category} of $\mathscr{A}$, where $\ast\in\{{\rm blank},-,+,b\}$.
The {\it left $\mathscr{X}$-derived category} $D_{\mathbb{L}\mathscr{X}}^{\ast}$($\mathscr{A})$ of $\mathscr{A}$ is defined dually.
\vspace{0.2cm}

{\bf Example 4.2.} Let $\mathscr{A}$ have enough projectives.

(1) If $\mathscr{X}=\mathscr{P(A)}$, then $D_{\mathbb{R}\mathscr{X}}^{*}$($\mathscr{A}$) is the usual derived category $D^{*}(\mathscr{A})$.

(2) If $\mathscr{X}=\mathscr{G(A)}$ (the subcategory of $\mathscr{A}$ consisting of
Gorenstein projective objects), then $D_{\mathbb{R}\mathscr{X}}^{*}$($\mathscr{A}$) is the Gorenstein
derived category $D_{gp}^{*}$($\mathscr{A}$) defined in [GZ].

\vspace{0.2cm}

The following two results are cited from [AHV].

\vspace{0.2cm}

{\bf Proposition 4.3.} {\it ([AHV])

(1) $D_{\mathbb{R}\mathscr{X}}^{-}(\mathscr{A})$ is a triangulated subcategory of
$D_{\mathbb{R}\mathscr{X}}(\mathscr{A})$, and $D_{\mathbb{R}\mathscr{X}}^{b}(\mathscr{A})$ is a triangulated subcategory
of $D_{\mathbb{R}\mathscr{X}}^{-}(\mathscr{A})$.

(2) For any $X^{\bullet}\in K^{-}(\mathscr{X})$ and $C^{\bullet}\in C(\mathscr{A})$, there exists an isomorphism of abelian groups:
$$\Hom_{K(\mathscr{A})}(X^{\bullet},C^{\bullet}) \cong  \Hom_{D_{\mathbb{R}\mathscr{X}}(\mathscr{A})}(X^{\bullet},C^{\bullet}).$$

(3) Let $\mathscr{X}\subseteq\mathscr{A}$ be admissible. Then the composition functor
$\mathscr{A}\to K^{b}(\mathscr{A})\to D_{\mathbb{R}\mathscr{X}}^{b}(\mathscr{A})$ is fully faithful,
where both functors are canonical ones.}

\vspace{0.2cm}

Set $$K^{-,\mathbb{R}\mathscr{X}b}(\mathscr{X}):=\{X^{\bullet}\in K^{-}(\mathscr{X})\mid\ {\rm there\ exists\ } n\in \mathbb{Z}\ {\rm such\ that}$$
$$\ \ \ \ \ \ \ \ \ \ \ \ \ \ \ \ \ \ \ \ \ \ \ \ H^i(\Hom_{\mathscr{A}}(X,X^{\bullet}))=0\ {\rm for\ any}\ X\in\mathscr{X}\ {\rm and}\ i \leq n\},$$ and
$$K^{+,\mathbb{L}\mathscr{X}b}(\mathscr{X}):=\{X^{\bullet}\in K^{+}(\mathscr{X})\mid\ {\rm there\ exists}\ n\in \mathbb{Z}\ {\rm such\ that}$$
$$\ \ \ \ \ \ \ \ \ \ \ \ \ \ \ \ \ \ \ \ \ \ \ \ H^i(\Hom_{\mathscr{A}}(X^{\bullet},X))=0\ {\rm for\ any}\ X\in\mathscr{X}\ {\rm and}\ i \leq n\}.$$

\vspace{0.2cm}

{\bf Proposition 4.4.} {\it ([AHV, Theorem 3.3]) If $\mathscr{X}$ is a contravariantly finite subcategory of $\mathscr{A}$,
then we have a triangle-equivalence $K^{-,\mathbb{R}\mathscr{X}b}(\mathscr{X})\simeq D_{\mathbb{R}\mathscr{X}}^{b}(\mathscr{A})$.}

\vspace{0.2cm}

As consequences of Propositions 4.4 and 4.3, we have the following two results.

\vspace{0.2cm}

{\bf Proposition 4.5.}  {\it Let ($\mathscr{X}$, $\mathscr{Y}$) be a balanced pair in $\mathscr{A}$. Then we have triangle-equivalences:
$$K^{-,\mathbb{R}\mathscr{X}b}(\mathscr{X})\simeq D_{\mathbb{R}\mathscr{X}}^{b}(\mathscr{A})= D_{\mathbb{L}\mathscr{Y}}^{b}(\mathscr{A})
\simeq K^{+,\mathbb{L}\mathscr{Y}b}(\mathscr{Y}).$$}

\vspace{0.2cm}

{\it Proof.} The first equivalence follows from Proposition 4.4, and the last one is its dual.

By [C1, Proposition 2.2], we have that the class of right $\mathscr{X}$-acyclic complexes coincides with that of left $\mathscr{Y}$-acyclic complexes
if ($\mathscr{X}$, $\mathscr{Y}$) is a balanced pair. Then $K^{b}_{\mathbb{R}\mathscr{X}{\text -}ac}(\mathscr{A})$ coincides with $K^{b}_{\mathbb{L}\mathscr{Y}{\text -}ac}(\mathscr{A})$,
and so we have $D_{\mathbb{R}\mathscr{X}}^{b}(\mathscr{A})= D_{\mathbb{L}\mathscr{Y}}^{b}(\mathscr{A})$. \hfill{$\square$}

\vspace{0.2cm}

For a balanced pair ($\mathscr{X}$, $\mathscr{Y}$) in $\mathscr{A}$, we call $D_{\mathbb{R}\mathscr{X}}^{b}(\mathscr{A})$ and
$D_{\mathbb{L}\mathscr{Y}}^{b}(\mathscr{A})$ the {\it relative bounded derived category}
relative to ($\mathscr{X}$, $\mathscr{Y}$), and denote them by $D_{\ast}^{b}(\mathscr{A})$.
The following result means that the relative cohomology group $\Ext_\ast^i(M,N)$ may be computed in the relative bounded derived category $D_{\ast}^{b}(\mathscr{A})$.

\vspace{0.2cm}

{\bf Proposition 4.6.} {\it Let ($\mathscr{X}$, $\mathscr{Y}$) be a balanced pair in $\mathscr{A}$.
Then for any $M,N\in\mathscr{A}$ and $i\geq1$, there exists an isomorphism of abelian groups:
$$\Ext_\ast^i(M,N)\cong \Hom_{D_{\ast}^{b}(\mathscr{A})}(M,N[i]).$$}

\vspace{0.2cm}

{\it Proof.} Let $\varepsilon:X_M^\bullet\to M$ be a $\mathscr{X}$-resolution of $M$. View $M$ as a stalk complex concentrated in degree zero.
Note that $\varepsilon$ is a right $\mathscr{X}$-quasi-isomorphism. So $M\cong X^\bullet_M$ in $D_{\mathbb{R}\mathscr{X}}(\mathscr{A})$.
Since $X^\bullet_M\in K^-(\mathscr{X})$, by Proposition 4.3 we have isomorphisms of abelian groups:
$$\Ext_\ast^i(M,N)=H^i\Hom_\mathscr{A}(X_M^\bullet,N)\cong\Hom_{K(\mathscr{A})}(X_M^\bullet,N[i])\cong
\Hom_{D_{\mathbb{R}\mathscr{X}}(\mathscr{A})}(X_M^\bullet,N[i])$$$$\cong\Hom_{D_{\mathbb{R}\mathscr{X}}(\mathscr{A})}(M,N[i])
\cong\Hom_{D_{\mathbb{R}\mathscr{X}}^b(\mathscr{A})}(M,N[i])\cong\Hom_{D_{\ast}^b(\mathscr{A})}(M,N[i]).$$ \hfill{$\square$}
\vspace{0.2cm}

Given a balanced pair ($\mathscr{X}$, $\mathscr{Y}$) in $\mathscr{A}$, from the viewpoint of relative derived category relative to ($\mathscr{X}$, $\mathscr{Y}$),
the $\mathscr{X}$-resolution of an object $M\in\mathscr{A}$ is exactly an isomorphism $X_M^\bullet\to M$ in $D_{\ast}^{b}(\mathscr{A})$,
where $X_M^\bullet\in K^-(\mathscr{X})$ with components vanish in the positive degrees, while the $\mathscr{Y}$-coresolution of an object
$N\in\mathscr{A}$ is exactly an isomorphism $N\to Y_N^\bullet$ in $D_{\ast}^{b}(\mathscr{A})$,
where $Y_N^\bullet\in K^+(\mathscr{Y})$ with components vanish in the negative degrees. In the following result,
we give some criteria for computing the $\mathscr{X}$-resolution dimension of an object in $\mathscr{A}$
in terms of the vanishing of relative cohomology groups.

\vspace{0.2cm}

{\bf Theorem 4.7.} {\it Let ($\mathscr{X}$, $\mathscr{Y}$) be an admissible balanced pair in $\mathscr{A}$.
Then the following statements are equivalent for any $M \in\mathscr{A}$ and $n\geq 0$.

(1) $\mathscr{X}$-$\resdim M\leq n$.

(2) $\Ext_{\ast}^{\geq n+1}(M,N)=0$ for any $N\in \mathscr{A}$.

(3) $\Ext_{\ast}^{n+1}(M,N)=0$ for any $N\in \mathscr{A}$.

(4) For any $\mathscr{X}$-resolution $X^\bullet\to M$ of $M$, we have $\Ker d_X^{-n+1}\in\mathscr{X}$.}

\vspace{0.2cm}

{\it Proof.} Both $(2)\Rightarrow (3)$ and $(4)\Rightarrow (1)$ are trivial.

$(1)\Rightarrow (2)$ Let $$0 \to X^{-n}\to X^{-n+1}\to \cdots\to X^0\to M\to 0$$ be an
$\mathscr{X}$-resolution of $M$. Then $\Hom_\mathscr{A}(X^{-i},N)$ =0 for any $N\in \mathscr{A}$ and $i\geq n+1$
and the assertion follows.

$(3)\Rightarrow (4)$ Let $$\cdots \to X^{-n}\buildrel {d_X^{-n}} \over \longrightarrow X^{-n+1}\to \cdots\to X^0\to M\to 0$$
be an $\mathscr{X}$-resolution of $M$. Then we have a $\ast$-acyclic sequence:
$$0\to \Ker d_X^{-n}\to X^{-n}\to \Ker d_X^{-n+1}\to 0. \eqno{(4.1)}$$ Since $\Ext_{\ast}^{n+1}(M,\Ker d_X^{-n})=0$, by the dimension shifting we have
$$\Ext_{\ast}^{1}(\Ker d_X^{-n+1}, \Ker d_X^{-n})\cong \Ext_{\ast}^{n+1}(M,\Ker d_X^{-n})=0.$$ It follows from [EJ2, Theorem 8.2.3] that (4.1) splits.
So $\Ker d_X^{-n+1}$ is isomorphic to a direct summand of $X^{-n}$ and $\Ker d^{-n+1}_X\in\mathscr{X}$. \hfill{$\square$}

\vspace{0.2cm}

Dually, we have the following

\vspace{0.2cm}

{\bf Theorem 4.8.} {\it Let ($\mathscr{X}$, $\mathscr{Y}$) be an admissible balanced pair in $\mathscr{A}$.
Then the following statements are equivalent for any $N \in\mathscr{A}$ and $n\geq 0$.

(1) $\mathscr{Y}$-$\coresdim N\leq n$.

(2) $\Ext_{\ast}^{\geq n+1}(M,N)=0$ for any $M\in \mathscr{A}$.

(3) $\Ext_{\ast}^{n+1}(M,N)=0$ for any $M\in \mathscr{A}$.

(4) For any $\mathscr{Y}$-coresolution $N\to Y^\bullet$ of $N$, we have $\Im d_Y^{n-1}\in\mathscr{Y}$.}

\vspace{0.2cm}

As an immediate consequence of Theorems 4.7 and 4.8, we have the following

\vspace{0.2cm}

{\bf Corollary 4.9.} ([C1, Corollary 2.5]) {\it Let ($\mathscr{X}$, $\mathscr{Y}$) be an admissible balanced pair in $\mathscr{A}$. Then
$\mathscr{X}$-$\resdim \mathscr{A}=\mathscr{Y}$-$\coresdim \mathscr{A}$.}

\vspace{0.2cm}

Set
$$\mathscr{X}{\text -}\resdim \mathscr{Y}:=\sup\{\mathscr{X}{\text-}\resdim Y\mid Y\in \mathscr{Y} \},$$ and
$$\mathscr{Y}{\text -}\coresdim \mathscr{X}:=\sup\{\mathscr{Y}{\text-}\coresdim X\mid X\in \mathscr{X} \}.$$

\vspace{0.2cm}

{\bf Theorem 4.10.} {\it For a balanced pair ($\mathscr{X}$, $\mathscr{Y}$) in $\mathscr{A}$, in $D_{\ast}^{b}(\mathscr{A})$ we have

(1) If $\mathscr{X}{\text -}\resdim \mathscr{Y}<\infty$, then $K^b(\mathscr{Y})\subseteq K^b(\mathscr{X})$.

(2) If $\mathscr{Y}{\text -}\coresdim \mathscr{X}<\infty$, then $K^b(\mathscr{X})\subseteq K^b(\mathscr{Y})$.}

\vspace{0.2cm}

{\it Proof.} (1) It suffices to show that for any $Y^\bullet\in K^b(\mathscr{Y})$, there exists a right $\mathscr{X}$-quasi-isomorphism
$X^\bullet_Y\to Y^\bullet$ with $X^\bullet_Y\in K^b(\mathscr{X})$. We proceed
by induction on the width $\omega(Y^\bullet)$ (:=the cardinal of the set $\{Y^i \neq 0\mid i\in \mathbb{Z}\}$) of $Y^\bullet$.

For the case $\omega(Y^\bullet)$=1, the assertion follows from the assumption that $\mathscr{X}$ is contravariantly finite and $\mathscr{X}$-$\resdim \mathscr{Y}<\infty$.

Let $\omega(Y^\bullet)\geq 2$ with $Y^j\neq 0$ and $Y^i=0$ for any $i<j$. Put $Y^\bullet_1$:=$Y^j[-j-1]$ and $Y^\bullet_2:=\sigma^{> j}Y^\bullet$.
Let $g=d_Y^j[-j-1]$ where $d_Y^j$ is the $j$th differential of $Y^\bullet$.
We have a distinguished triangle $$Y^\bullet_1\buildrel {g} \over \longrightarrow Y^\bullet_2\to Y^\bullet\to Y^\bullet_1[1]$$ in $K^b(\mathscr{Y})$.
By the induction hypothesis, there exist right $\mathscr{X}$-quasi-isomorphisms
$f_{Y_1}$: $X_{Y_1}^\bullet\to Y^\bullet_1$ and $f_{Y_2}$: $X_{Y_2}^\bullet\to Y^\bullet_2$ with $X_{Y_1}^\bullet$, $X_{Y_2}^\bullet$ $\in K^{b}(\mathscr{X})$.
Then by Lemma 2.4, $f_{Y_2}$ induces an isomorphism:
$$\Hom_{K^b(\mathscr{A})}(X_{Y_1}^\bullet,X_{Y_2}^\bullet)\cong \Hom_{K^b(\mathscr{A})}(X_{Y_1}^\bullet,Y^\bullet_2).$$
So there exists a morphism $f:X_{Y_1}^\bullet \to X_{Y_2}^\bullet$, which is unique up to homotopy, such that
$f_{Y_2} f=gf_{Y_1}$. Put $X_{Y}^\bullet=\Con(f)$. We have the following distinguished triangle
\begin{center}
$X_{Y_1}^\bullet \buildrel {f} \over \longrightarrow X_{Y_2}^\bullet\to X_{Y}^\bullet\to X_{Y_1}^\bullet[1].$
\end{center}
in $K^b(\mathscr{X})$. Then there exists a morphism $f_{Y}:X_{Y}^\bullet\to Y^\bullet$ such that the following diagram commutes:
$$\xymatrix{X_{Y_1}^\bullet \ar[r]^{f}\ar[d]^{f_{Y_1}} & X_{Y_2}^\bullet \ar[r] \ar[d]^{f_{Y_2}} &
X_{Y}^\bullet \ar[r]\ar@{-->}[d]^{f_{Y}} & X_{Y_1}^\bullet[1] \ar[d]^{f_{Y_1}[1]} \\
Y^\bullet_1 \ar[r]^{g} & Y^\bullet_2 \ar[r] &
Y^\bullet \ar[r] & {Y^\bullet_1}[1].}$$
For any $X\in \mathscr{X}$ and $n\in \mathbb{Z}$, we have the following commutative diagram with exact rows:
$$\xymatrix{(X,X_{Y_1}^\bullet) \ar[r]\ar[d]^{(X, f_{Y_1})} & (X,X_{Y_2}^\bullet) \ar[r] \ar[d]^{(X, f_{Y_2})} &
(X,X_{Y}^\bullet) \ar[r]\ar@{-->}[d]^{(X, f_{Y})} & (X,X_{Y_1}^\bullet[1]) \ar[d]^{(X, f_{Y_1}[1])}\ar[r] & (X,X_{Y_2}^\bullet[1]) \ar[d]^{(X,f_{Y_2}[1])}\\
(X,Y^\bullet_1) \ar[r] & (X,Y^\bullet_2) \ar[r] &
(X,Y^\bullet) \ar[r] & (X,{Y^\bullet_1}[1])\ar[r] & (X,{Y^\bullet_2}[1]),}$$
where $(X,-)$ denotes the functor $\Hom_{K(\mathscr{A})}(X,[n](-))$. Since $f_{Y_1}$ and $f_{Y_2}$ are right $\mathscr{X}$-quasi-isomorphisms,
we have that $(X,f_{Y_1})$ and $(X,f_{Y_2})$ are isomorphisms. So $(X,f_{Y})$ is also an isomorphism and $f_{Y}$ is a right $\mathscr{X}$-quasi-isomorphism.
The proof is finished.

(2) It is dual to (1). \hfill{$\square$}

\vspace{0.2cm}

Let $A$ be a finite-dimensional algebra over a field $k$. We use $\mod A$ to denote the category of finitely generated left $A$-modules,
and use $\proj A$ (resp. $\inj A$) to denote the full subcategory of $\mod A$ consisting of projective (resp. injective) modules.
For a module $M\in\mod A$, we use $\pd_AM$ and $\id_AM$ to denote the projective and injective dimensions of $M$ respectively.
As an application of Theorem 4.10, we get the following

\vspace{0.2cm}

{\bf Corollary 4.11.} ([H]) {\it For a finite-dimensional algebra over a field $k$, in $D^b(A)$ we have

(1) $\pd_A D(A_A)<\infty$ if and only if $K^b(\inj A)\subseteq K^b(\proj A)$.

(2) $\id _AA<\infty$ if and only if $K^b(\proj A)\subseteq K^b(\inj A)$.

(3) $A$ is Gorenstein if and only if $K^b(\proj A)= K^b(\inj A)$.}

\vspace{0.2cm}

{\it Proof.} We only prove (1), because (2) is dual to (1), and (3) is an immediate consequence of (1) and (2).

The necessity follows from Theorem 4.10. For the sufficiency, since $K^b(\inj A)\subseteq K^b(\proj A)$ in $D^b(A)$,
we have $D(A_A)\in K^b(\proj A)$. Then there exists a quasi-isomorphism $Q^\bullet\to D(A_A)$ with $Q^\bullet\in K^b(\proj A)$.
Let $P^\bullet\to D(A_A)$ be the projective resolution of $D(A_A)$ in $\mod A$. It follows that $P^\bullet$ and $Q^\bullet$ are homotopy equivalence.
Thus $P^\bullet\in K^b(\proj A)$ and hence $\pd_A D(A_A)<\infty$.  \hfill{$\square$}

\vspace{0.2cm}

Let ($\mathscr{X}$, $\mathscr{Y}$) be a balanced pair in $\mathscr{A}$. It follows from Proposition 4.3 that $K^b(\mathscr{X})$ is a triangulated subcategory of
$D_{\ast}^{b}(\mathscr{A})$. Motivated by the definition of classical singularity categories, we introduce the following

\vspace{0.2cm}

{\bf Definition 4.12.} We call the quotient category $D_{\mathbb{R}\mathscr{X}{\text-}sg}(\mathscr{A}):= D^{b}_\ast(\mathscr{A})$ /$K^{b}(\mathscr{X})$
the {\it right $\mathscr{X}$-singularity category} relative to ($\mathscr{X}$, $\mathscr{Y}$), and call $D_{\mathbb{L}\mathscr{Y}{\text-}sg}(\mathscr{A}):
=D^{b}_\ast(\mathscr{A})$ /$K^{b}(\mathscr{Y})$ the {\it left $\mathscr{Y}$-singularity category} relative to ($\mathscr{X}$, $\mathscr{Y}$).

\vspace{0.2cm}

Let $A$ be a finite-dimensional algebra over a field $k$. In the case for $\mathscr{X}=\proj A$, we have that $D^{b}_\mathscr{X}(\mathscr{A})$ coincides with the usual
bounded derived category $D^{b}(\mathscr{A})$ and $D_{\mathbb{R}\mathscr{X}{\text-}sg}(\mathscr{A})$ is the classical singularity category
$D_{sg}(A)$ which is called the ``stabilized derived category" in [Bu]. For the properties of singularity categories and related topics,
we refer to [C2], [CZ], [H], [O], [R], and so on. It is known that $D_{sg}(A)=0$ if and only if $A$ is of finite global dimension. So $D_{sg}(A)$ measures the homological
singularity of the algebra $A$.

\vspace{0.2cm}

{\bf Theorem 4.13.} {\it Let ($\mathscr{X}$, $\mathscr{Y}$) be a balanced pair in $\mathscr{A}$.

(1) If $\mathscr{X}$-$\resdim \mathscr{Y}<\infty$
and $\mathscr{Y}$-$\coresdim \mathscr{X}<\infty$, then $D_{\mathbb{R}\mathscr{X}{\text-}sg}(\mathscr{A})=D_{\mathbb{L}\mathscr{Y}{\text-}sg}(\mathscr{A})$.

(2) If ($\mathscr{X}$, $\mathscr{Y}$) is admissible and $\mathscr{X}$-$\resdim \mathscr{A}< \infty$,
then $D_{\mathbb{R}\mathscr{X}{\text-}sg}(\mathscr{A})=0=D_{\mathbb{L}\mathscr{Y}{\text-}sg}(\mathscr{A})$.}

\vspace{0.2cm}

{\it Proof.} (1) If $\mathscr{X}$-$\resdim \mathscr{Y}<\infty$ and $\mathscr{Y}$-$\coresdim \mathscr{X}<\infty$, then it follows from Theorem 4.10 that
$K^b(\mathscr{X})= K^b(\mathscr{Y})$. So we have $D_{\mathbb{R}\mathscr{X}{\text-}sg}(\mathscr{A})=D_{\mathbb{L}\mathscr{Y}{\text-}sg}(\mathscr{A})$.

(2) Since ($\mathscr{X}$, $\mathscr{Y}$) is an admissible balanced pair and $\mathscr{X}$-res.dim $\mathscr{A}< \infty$,
it follows from Corollary 4.9 that $\mathscr{Y}$-$\coresdim \mathscr{A}<\infty$. For the first equality it suffices to show that for any
$A^\bullet\in K^b(\mathscr{A})$, there exists a right $\mathscr{X}$-quasi-isomorphism $X^\bullet_A\to A^\bullet$ with $X^\bullet_A\in K^b(\mathscr{X})$.
By using an induction on the width $\omega(A^\bullet)$ of $A^\bullet$ and a similar argument to that in proof of Theorem 4.10, we get the assertion.
Dually, we get the second equality.   \hfill{$\square$}

\vspace{0.2cm}

Let $A$ be a finite-dimensional algebra over a field $k$. We use $\Gproj A$ (resp. $\Ginj A$) to
denote the full subcategory of $\mod A$ consisting of Gorenstein projective
(resp. injective) modules. It follows from [C1, Proposition 2.6] that ($\Gproj A$, $\Ginj A$) is an
admissible balanced pair in $\mod A$ whenever $A$ is Gorenstein. Let ($\mathscr{X}$, $\mathscr{Y}$)=($\Gproj A$, $\Ginj A$).
Consider the following quotient categories (cf. [GZ]):
$$D_{\mathbb{R}\mathscr{G}}{\text -}sg(A):=D_{*}^b(\mod A)/K^b(\Gproj A),$$
$$D_{\mathbb{L}\mathscr{G}}{\text -}sg(A):=D_{*}^b(\mod A)/K^b(\Ginj A).$$

By Theorem 4.13(2), we have the following

\vspace{0.2cm}

{\bf Corollary 4.14.} {\it Let $A$ be a Gorenstein algebra. Then $D_{\mathbb{R}\mathscr{G}}{\text -}sg(A)=0$ and
$D_{\mathbb{L}\mathscr{G}}{\text -}sg(A)=0$.}

\vspace{0.5cm}

{\bf Acknowledgements.} This research was partially supported by NSFC (Grant No. 11171142) and a Project
Funded by the Priority Academic Program Development of Jiangsu Higher Education Institutions.

\end{document}